\PassOptionsToPackage{pdfpagelabels=false}{hyperref} 
\documentclass{ws-rv9x6}
\usepackage[OT4,T1]{fontenc}
\usepackage[utf8]{inputenc}
\usepackage{subfigure}   % required only when side-by-side / subfigures are used
\usepackage{ws-rv-thm}   % comment this line when `amsthm / theorem / ntheorem` package is used
\usepackage{ws-rv-van}   % numbered citation & references (default)
\usepackage{float}
\RequirePackage[colorlinks,citecolor=blue,urlcolor=blue]{hyperref}
\usepackage{graphicx}
\usepackage{epsfig}
\usepackage{epstopdf}
  \hypersetup{
    pdftitle={Full vs. no ...}, %%<--To wymienić
    pdfauthor={Marek Skarupski and Krzysztof Szajowski}, %%<--TO wymienić
    colorlinks,
    urlcolor=blue,
    filecolor=magenta,
    citecolor=green,
    linkbordercolor={1 1 1}, % set to white
    citebordercolor={1 1 1},  % set to white
    urlbordercolor={ 1 1 1}  % set to white
  }

\newcommand{\one}{\mathbb{I}}
\makeindex
\begin{document}

\chapter[Full vs. no information BCG with FI]{\label{ra_ch1} Full vs. no information best choice game\\
 with finite horizon}

\author[M. Skarupski and K. Szajowski]{Marek Skarupski\footnote{Marek Skarupski: Wroc\l{}aw University of Science and Technology, Faculty of Pure and Applied Mathematics,  Wybrze\.ze Wyspia\'nskiego 27, PL-50-370 Wroc\l{}aw.\\ \emph{E-mail: } Marek.Skarupski@pwr.edu.pl} and Krzysztof Szajowski\footnote{Krzysztof Szajowski}}
%\index[aindx]{Author, F.} % or \aindx{Author, F.}
%\index[aindx]{Author, S.} % or \aindx{Author, S.}

\address{Wroc\l{}aw University of Science and Technology\\
Krzysztof.Szajowski@pwr.edu.pl\footnote{Wroc\l{}aw University of Science and Technology, Faculty of Pure and Applied Mathematics,   Wybrze\.ze Wyspia\'nskiego 27, PL-50-327 Wroc\l{}aw.}}

\begin{abstract}
Let us consider a two companies A and B. Both of them are interested in buying a set of some goods. The company A is a big corporation and it knows the actual value of the good on the market and is able to observe the previous values of them. The company B has no information about the actual value of the good but it can compare the actual position of the good on the market with the previous position of the good offered. Both of the players want to choose the very best object overall. The recall is not allowed. The number of the objects is fixed and finite. One can think about these two types of buyers a business customer vs. an individual customer. The mathematical model of the competition between them is presented and the solution is defined and constructed.
\end{abstract}

%\markright{Customized Running Head for Odd Page} % default is Chapter Title.
\body

%\tableofcontents
\section{Introduction.}
The very well known secretary problem has also many modifications. Ferguson\cite{Ferguson1989Who} has made a review of the concepts of the best choice problem going back to the age of Kepler and Cayley. Presman and Sonin\cite{PresmanSonin72Thebest} considered so called no-information problem in which the appearing objects comes from the rank distribution, i.e. that the objects are observable, the decision maker can rank them and all permutations of the appearing objects are equally possible. Another approach was presented by Gilbert and Mosteller\cite{GilbertMosteller66Recognizing} where the exact value of the object is observable, and the distribution of the object is known (it is assumed to be an uniform distribution on the interval $[0,1]$). Both of the ideas can be described as an optimal stopping of Markov chain. In both of them there is only one decision maker and there is no competition concept. The game concept of the secretary problem was introduced by Dynkin\cite{Dyn69:Game}. Many examples were solved by Yasuda\cite{YasudaNakagamiKurano82,Yasuda85:Neveu}. 

\subsection{Business motivation.}
Consider a two companies A and B. Both of them are interested in buying a set of some goods, ex. an asset on a stock exchange. Th company A is a big corporation and it knows the actual value of the good on the market. What is more it knows the previous values of the objects and can compare them together. The problem of the company B is that it has no information about the actual value of the good. However the owner of the company B can compare the actual position of the good in the market with the previous observations. Both of the players want to choose the very best object overall without possibility of recall. The number of the objects is fixed and finite. A very good example can be described from the position of reliability. Consider a two buyers of the same item. Both of them want to buy the most reliable object. The buyer A has possibility to get know the reliability function values derived by the experts and quality controllers. The player B has no such a contact and intelligence, so he must rely on his basic knowledge and the knowledge of the previous observation, i.e. he can judge is the object better or worse than the previous one. We can say that the buyers of the objects are two types: first is a business customer and second is an individual customer.

\subsection{Related game models.} The considered game recalls various conflict driven by the stochastic sequences. Usually the bilateral setting of the decision problem is preceded by unilateral consideration. The form of the optimal strategy in the decision problem is the inspiration for mean-value formulation. The threshold strategies are crucial tool for the optimal stopping problems. The simplest case related to observation of a sequences of random variables can be found e.g. in papers [\refcite{Sza84:maxima}] or [\refcite{Por90:overMAXIMA}]. The bilateral extension of these models can be found in the papers [\refcite{NeuPorSza96:imperfect}]. Two players, $I$ and $II$, observe sequentially a known finite number (or a number having a geometric distribution) of independent and identically distributed random variables. They must choose the largest. The variables cannot be perfectly observed. When a random variable is sampled, the sampler is informed only whether it is greater or less than some level he has specified. Each player can choose at most one observation. After the sampling, the players decide for acceptance or rejection of the observation. If both accept the same observation, Player 1 has the priority. The class of adequate strategies and a gain function are constructed. In the finite case, the game has a solution in pure strategies. In the case of a geometric distribution, Player 1 has a pure equilibrium strategy and Player 2 has either a pure equilibrium strategy or a mixture of two pure strategies. The game is symmetrical as the players are watching the same string to the same extent.  Increasing the opposing interests is possible by the completely different preferences of the players. Evaluation of the same object by two decision makers can mean that players observe the different coordinates of the vector and formulate their expectations for their realization. When players aim is to achieve a minimum level of the observed rate then the problem can be reduced to a game in which strategies are setting of just levels. Discussion of such issues can be found in the works of Sakaguchi (e.g. [\refcite{Sak73:bivariateJORSJ}]). However and in those tasks, though, that the information players are incomplete, lacking clear asymmetry players. The pay-offs of the players are function of the thresholds and the perfect comparison of the observed variable with these defined levels is guarantee. Asymmetric tools in measure of the observed r.v. are presented by Sakaguchi and the Second Author\cite{SakSza00:bivariate}. However, for private random variables. These players with asymmetric tools applied to there same sequence are subject of consideration in the paper.   
%Gra jest symetryczna jako że gracze obserwują ten sam ciąg w takim samym zakresie. Zwiękaszenie przeciwstawnych interesów jest możliwe przez odmienne preferencje graczy. Ocena tego samego obiektu może oznaczać, że gracze obserwują różne wektory i formułują swoje oczekiwania względem ich realizacji. Gdy celem graczy jest osiągnięcie minimalnego poziomu obserwowanego wskaźnika zadanie można sprowadzić do gry w której strategiami są te własnie poziomy. Omówienie takich zagadnień można znaleźć w pracach Sakaguchi. Jednakże i w tych zadaniach, mimo, że informacje graczy są niepełne, brakuje wyraźnej asymetrii graczy.
 
\subsection{Mathematical formulation of the problem.} In fine-tuning the mathematical model, we will use the methods of optimal stopping of stochastic processes \cite{choRobSig71:Expectation} for Markov sequences\cite{DynkinYushkievich67,Shi78:OSR} and game models with optimal stopping\cite{Dyn69:Game} of such sequences, similar to what is done in the works [\refcite{Sza94:ZOR,Sza95:SIAM}].

Let $(\Omega,\mathcal{F},\mathbf{P})$ be enough reach probability space to define the random sequence $\{X_n\}_{n=0}^N$, $X_\cdot:\Omega\rightarrow\mathbb{R}\subset\Re$, $N\in \mathbb{N}\cup\{\infty\}$. There are two observers (and at the same time decision makers) of the basic sequence defined by the mappings $\{\varphi^i_n\}$, $i=1,2$, where $\varphi^i_n:\Re^n\rightarrow\Re$, having his aims defined by the pay-off functions $f^i:\mathbb{R}\rightarrow\Re$. Other words, the player $i$ at moment $n$ observes $\xi^i_n=\varphi^i_n(X_1,\ldots,X_n)$. The strategies of the players are stopping times $\tau\in\mathfrak{S}^i$ with respect to the appropriate filtrations $\mathcal{F}^i_n=\sigma\{\xi^i_1,\ldots,\xi^i_n\}$. Each player, based on the observations available to him, is tasked with choosing the moment of accepting the state of the process based on the previous observations so as to maximize the expected payment. 
\begin{equation}\label{expPayOff}
\hat{v}^i=\sup_{\tau^i\in\mathfrak{S}^i}\mathbb{E}f^i(X_{\tau^i}).
\end{equation}
It can be used to reduce the initial problem to the task of optimal retention of conditional expected values relative to its filtration\cite{Shi78:OSR}. Let us calculate for every $n\in\mathbb{N}$
\begin{subequations}

\begin{equation*}\label{ConExpPayOff1}
\hat{f}^i(\vec{\xi}^i_n)=\mathbf{E}[f^i(X^i_n)|\mathcal{F}^i_n].
\end{equation*}
We have
\begin{equation*}\label{ConExpPayOff2}
\hat{v}^i=\sup_{\tau^i\in\mathfrak{S}^i}\mathbf{E}\hat{f}^i(\vec{\xi}^i_{\tau^i}).
\end{equation*}
\end{subequations}
Let us assume that the observation processes $\{\xi^i_n\}_{n=0}^N$, $i=1,2$, belongs to Markov processes. In this case the solution of the problem \eqref{expPayOff} can be get by the procedure described by Shiryaev\citet[Ch.3]{Shi78:OSR} which is based on the Bellman-Jacobi equation. Denote $\mathfrak{S}_n=\{\tau\in\mathfrak{S}:\tau\geq n\}$. When we have two decision makers hunting for a convenient state of the process have the right to declare the stopping at most twice. The second only, when the first hired state is assigned to the opponent. The natural set of strategies are $\mathfrak{U}^i=\{(\tau^i,\{\sigma^i_n\}_{n=0}^N):\tau^i\in\mathfrak{S},\sigma^i_n\in\mathfrak{S}_n^i\}$. The pay-off in the competitive case is defined in various way. Following the discussion in [\refcite{Sza95:SIAM}] for given $\rho^i\in\mathfrak{U}^i$, $i=1,2$,
\begin{align}\label{GamePayOff}
K_1(\rho^1,\rho^2)=\mathbf{E}&\Big[\one_{\{\tau^1<\tau^2\}}\big(\hat{f}^1(\xi^1_{\tau^1})-\hat{v}^2(\tau^1,\xi^2_{\sigma^2_{\tau^1}})\big)\\
\nonumber&\quad +\one_{\{\tau^1=\tau^2\}}\big[p\big(\hat{f}^1(\xi^1_{\tau^1})-\hat{v}^2(\tau^1,\xi^2_{\sigma^2_{\tau^1}})\big)\\
\nonumber&\quad+(1-p)\big(\hat{v}^1(\tau^2,\xi^1_{\sigma^1_{\tau^2}})-\hat{f}^2(\xi^2_{\tau^2})\big)\big]\\
\nonumber&\quad+\one_{\{\tau^1>\tau^2\}}\big(\hat{v}^1(\tau^2,\xi^1_{\sigma^1_{\tau^2}})-\hat{f}^2(\xi^2_{\tau^2})\big)\Big],
\end{align}
where $\one_A$ is the characteristic function of $A$ and 
\begin{equation*}
\hat{v}^i(n,\xi^i_n)=\sup_{\tau^i\in\mathfrak{S}^i_n}\mathbf{E}\hat{f}^i(\vec{\xi}^i_{\tau^i})
\end{equation*}
and $0\leq p \leq 1$ is the priority parameter, i.e. the probability that the state will be assigned to the player 1. The pair of strategies $({\rho^1}^\star,{\rho^2}^\star)$ is the solution of the problem if for every $\rho^i\in\mathfrak{U}^i$
\[
K_1({\rho^1}^\star,{\rho^2}^\star)\geq K_1(\rho^1,{\rho^2}^\star) \text{ and }
K_1({\rho^1}^\star,{\rho^2}^\star)\geq K_1({\rho^1}^\star,\rho^2).
\]
In practice it is difficult to construct the solution and calculate the value of the problem in such general form . However, for some natural cases each player can estimate his final reward by calculating his potential reward (award) based on his knowledge (filtration). The idea of these simplifications is presented in the next sections.
  
\section{Formulation of the game} % Initial capital letter, then lower case. No full stop.
\subsection{The description of the model}
Consider a game in which two players want to choose the best object overall. They observe $N$ objects sequentially. They get a profit only if the player chooses the best object and the rival will not find the better one. In other case he gets the award. If both players wants to stop on the current object the nature chooses it by the fair coin toss.
Suppose that:
\begin{enumerate}
\item The player I has no information, i.e. he observes only the relative ranks of the current objects.
\item The player II has full information, i.e. he observes sequentially $X_1,...,X_N$ i.i.d., sees its value, and also can calculate the rank of the current object.
\end{enumerate} 
To be more specific let us denote by $Y_n$ the relative rank of the $n$-th observation
\begin{equation}
Y_n = \# \lbrace 1 \leq i \leq n : X_i \leq X_n \rbrace .
\end{equation}
The filtration of the player II is $\mathcal{F}^{(1)}_n=\sigma(X_1,...,X_n)$ and for the player I is $\mathcal{F}^{(2)}_n=\sigma(Y_1,...,Y_n)$. Note that $\mathcal{F}^{(2)}_n \subset \mathcal{F}^{(1)}_n$ for every $n$. Denote by $\mathcal{T}_1$ a set of all stopping times with respect to the family $\lbrace \mathcal{F}_{n} \rbrace_{n=1}^{N}$. Let $\mathcal{T}_{1}^{0}$ denote a set of all stopping times $\tau \in \mathcal{T}_1$ such that
$X_{n} = \max \lbrace X_{1},...,X_{n} \rbrace$  on  $\lbrace \tau = n \rbrace$, $n=1,\ldots,N$.
Define the moments where the biggest observation appear, i.e. $\tau_{1}=1$,
$\tau_{k}= \inf \lbrace n: \tau_{k-1} \leq n \leq N, X_{n} = \max \lbrace X_{1},...,X_{n} \rbrace \rbrace$ for $k=1,\ldots,N$. We observe that the sequence $\tau_1, \tau_2,... \in \mathcal{T}_{1}^{0}$.
Now let us consider the following chain
\begin{displaymath}
Z_{k}= (\tau_k, X_{\tau_k}) \textrm{ on } \lbrace \tau_k < N+1 \rbrace, Z_{k}=(\tau_{N+1},\partial), 
\end{displaymath}
where $\partial$ is special absorbing state. It is easy to see that $\lbrace Z_{k} \rbrace_{k=1}^{N+1}$ is a Markov chain with transition probabilities (cf. [\refcite{Bojdecki78Onoptimal}])
\begin{equation}
p((n,x),(m,B)) = %P(\tau_{k+1}=m, X_{m} \in B | \tau_{k}=n, X_n=x)= 
x^{m-n-1}\int_{B}dy, m>n
\end{equation} 
and $0$ otherwise, with $B \subseteq (x,1]$.

The reward for the player II for stopping on $n$th object of the value $X_n=x$ is 
\begin{equation}
s_{2,n}(x)=x^{N-n}
\end{equation}
and for continuing observation is given by Gilbert \& Mosteller\cite{GilbertMosteller66Recognizing,Bojdecki78Onoptimal} %(cf. Bojdecki\cite{Bojdecki78Onoptimal})
\begin{equation*}
%\begin{split}
c_{2,n}(x)=\sum_{k=n+1}^{N} p((n,x),(k,(x,1])) s_{1,k}(y) %=\sum_{k=n+1}^{N} x^{k-n-1} \int_{x}^{1}y^{N-k}dy  \\
= \sum_{k=n+1}^{N}\dfrac{x^{k-n-1}(1-x^{N-n})}{N-k+1}.
%\end{split}
\end{equation*}
In similar way for the player I consider a sequence of indicators  $\lbrace I_{n} \rbrace _{n=1}^{N}$, where $I_{k}=\one_{\{Y_k=1\}}$.
Let us denote by $\mathcal{G}_{n} = \sigma(I_{1},...I_{n})$ sequence of sigma fields generated by indicators and let $\mathcal{T}_{2}$ be the set of all stopping moments $\tau$ with respect to $\sigma$ -fields $\mathcal{G}_{n}, n=1,...,N$.
Define a process $\xi_t$ in the following way
$$\xi_{t}=\inf \lbrace n\geq \xi_{t-1}: I_{n}=1 \rbrace$$
with initial point $\xi_{0}=1$. Calculate transition probabilities\cite{DynkinYushkievich67}
\begin{equation}\label{equation2}
%\begin{split}
p_{n,m}  = P(\xi_{k+1}=m|\xi_{k}=n) %=\dfrac{P(\xi_{k+1}=m,\xi_{k}=n)}{P(\xi_{k}=n)} = \\
%& = \dfrac{P(I_{n}=1,I_{n+1}=...=I_{m-1}=0, I_{m}=1)}{P(I_{n}=1)} = \dfrac{1}{m} \prod_{k=n+1}^{m-1}\dfrac{k-1}{k} \\
 = \dfrac{n}{m(m-1)}.
%\end{split}
\end{equation}
The first player's reward for stopping on $n$th candidate (i.e. $Y_n=1$) is $s_{1,n}=\dfrac{n}{N}$ and for continuing observations
\begin{equation*}
c_{1,n}=\sum_{k=n+1}^{N}\dfrac{n}{k(k-1)}\dfrac{k}{N} = \dfrac{n}{N}\sum_{k=n+1}^{N}\dfrac{1}{k-1}.
\end{equation*}

\subsection{Equilibrium states.}
Suppose that we are in some moment $n$ and the value of the current candidate is $x$ and both players want to stop. If the player II gets the object (with probability $1-p$) he gets reward $s_{1,n}(x)$. With probability $p$ the player I gets the object so II must continue the observations and gets reward $c_{1,n}(x)$. The situation that in the future the opponent will find the best object is also included into the reward. Similar consideration gives us the payoff for the player I. Let us denote:
\begin{equation}
w_{2,n}(x) = \left( x^{N-n}-x^{N-n}\sum_{j=n+1}^{N}\dfrac{x^{j-N-1}-1}{N-j+1} \right)
\end{equation}
and
\begin{equation}
w_{1,n}=\dfrac{n}{N}\left(1-\sum_{j=n+1}^{N}\dfrac{1}{j-1}\right).
\end{equation}
Then the payoff matrix is given by
\begin{equation} (v_{1,n},v_{2,n}(x))=
\begin{array}{c|c|c}
  I \backslash II & \text{S} & \text{F} \\
  \hline 
 \text{S} &(2p-1)w_{1,n};(1-2p)w_{2,n}(x) & w_{1,n};-w_{2,n}(x)\\
  \hline
 \text{F} & -w_{1,n};w_{2,n}(x) & Tv_{1,n},Tv_{2,n}(x) \\
\end{array} 
\end{equation}
where $T$ stands for the Bellman operator\cite{Shi78:OSR}. %
Since both players want to maximize their profits we have the following conditions) to (S,S) be a Nash equilibrium:
\begin{subequations}
\begin{equation}\label{SS1}
\begin{cases}
(2p-1)w_{1,n}\geq -w_{1,n}\\
(1-2p)w_{2,n}(x)\geq w_{2,n}(x)
\end{cases}
\end{equation}
which leads to the equations
\begin{equation}\label{SS2}
\begin{cases}
\sum_{j=n+1}^{N}\dfrac{x^{j-N-1}-1}{N-j+1}\leq 1\\
\sum_{j=n+1}^{N}\dfrac{1}{j-1}\leq 1
\end{cases}
\end{equation}
\end{subequations}
For player I we get the optimal moments for stop. They are those number $n\geq n^*$ where $n^*$ is given by the standard optimal threshold (cf. [\refcite{GilbertMosteller66Recognizing}]):
\begin{equation}
n^* = \max \{ 0\leq n \leq N:  \sum_{k=n+1}^{N}\dfrac{1}{k-1} > 1 \}.
\end{equation}
and the optimal stopping set for player II is also the same set form the standard optimal stopping problem
\begin{equation}
D = \{ (n,x) \in \{1,2,...,N\}\times[0,1] : X_n=x, x\geq x_n  \}
\end{equation}
where $x_n$ is the solution of equation 
\begin{equation}
\begin{split}
\sum_{k=1}^{N-n}\dfrac{x^{-k}-1}{k} & =1.
\end{split}
\end{equation}
The optimal stopping times are:\\
for the player I: $ \tau_1 = \inf\lbrace n>n^*, Y_n=1  \rbrace$ and for the player II:  $\tau_2 = \inf \lbrace n: (n,X_n): X_n=\max\lbrace X_1,...,X_n  \rbrace=x \geq x_n  \rbrace. $.\\
\begin{lemma}
In the game described above strategy $(S,S)$ is pure Nash equilibrium if $X_n$ is a local maximum and $X_n\geq x_n, n> n^*$.
\end{lemma}

Let us consider case when $p<0.5$. Suppose that $n>n^*$ and the value of the current observation is $X_m = x \leq x_n$ and its relative rank is $1$. If we are bellow this threshold then it is more optimal for player II to change his strategy on F. The best response for player I on the strategy of the opponent is continue stopping if the expected future reward is not greater than the actual reward:
$$w_{1,n} \geq Tv_{1,n} .$$
To calculate the expected future reward we have to calculate the Bellman operator for payoff $Tv_{1,n}$. The no-information player knows that the opponent has more information. Since the opponent chooses the strategy F we know that the present value of the object is less than the threshold $x_n$. Suppose for a moment that we know this value and it is $x$. If we knew this value the future payoff would be
\begin{equation}
T(v_{1,n}|x) = \sum_{k=n+1}^{N}\left( x^{k-n-1} \int_{x}^{x \vee x_k} w_{1,k} dy + \int_{x\vee x_k}^{1} (2p-1)w_{1,k} dy\right).
\end{equation}
where $a\vee b = \max \{a,b \}$. However we have to average it. Knowing that the actual value is uniformly distributed on the interval $[0,x_n]$ (since the opponent wishes to continue the observations) we get
\begin{equation}
Tv_{1,n} = \dfrac{1}{x_n}\int_{0}^{x_n} T(v_{1,n}|x) dx.
\end{equation}
Let us consider the set $ M_1=\{ n^*<n\leq N: Tv_{1,n} \leq w_{1,n} \}$.
Note that this set in not empty. It contains the number $\{ N\}$. Using method of backward induction we can find the lower bound for this set, i.e. the index $ \tilde{n}= \min \{n^*\leq n\leq N: Tv_{1,n} \leq w_{1,n}\}$.
\begin{lemma}
Suppose that the current state of the process $(n,X_n)$ is such that $n\geq \tilde{n}, X_n=x \leq x_n$. and $X_n$ is a local maximum. Then the strategy $(S,F)$ is the pure Nash equilibrium in the game described above.
\end{lemma}

Now suppose that  $n=\tilde{n}-1$ and the current state of the process is $(n,X_n=x)$ where  $x\leq x_n$. Since the player I changes his strategy into $F$ it is necessary to check whether the condition $Tv_{2,n}(x) \geq w_{2,n}(x)$ is satisfied. Indeed it is true. Bellman's operator is the expected value of the future reward. Since now the reward $w_{2,n}(x)<0$ and the future reward is positive for $p<0.5$ it is more optimal to take an action $F$ for the player II. Now the same consideration are made for $\tilde{n}-2, \tilde{n}-3,...$ etc. So from this considerations we have the following
\begin{lemma}
Suppose that the current state of the process $(n,X_n)$ is such that $n\leq \tilde{n}, X_n=x \leq x_n$ and $X_n$ is a local maximum. Then the strategy $(F,F)$ is the pure Nash equilibrium in this state in the game described above.
\end{lemma}

Now consider the case when $n=n^*-1$ and $X_n=x > x_n$ is the local maximum. This is the opposite situation when the player II prefers to stop but the player I prefers to continue the observations. To find is strategy $(F,S)$ the equilibrium point we have to check is the condition
$$Tv_{2,n}(x)\leq w_{2,n}(x). $$
This is equivalent to the condition
\begin{equation}
2p\sum_{k=1}^{N-n}\dfrac{x^{-k}-1}{k} \geq -(1-2p)\sum_{k=1}^{N-n}\sum_{j=1}^{k-1}\left( \dfrac{j}{k}x^{-k} - \dfrac{x^{-j}}{k-j}+\dfrac{1}{k} \right)
\end{equation}
The expression under the double sum is positive in the interval $[0,1]$. It means that on the left hand side we have a positive number which is always bigger than the expression on the right hand side which is negative. There fore for $n=n^*-1$ and $x>x_n$ it is better for player II to not change his strategy. Continuing this calculations we get that it is also better to not change his strategy when $n< n^*$ and $x>x_n$.
\begin{lemma}
Suppose that the current state of the process $(n,X_n)$ is such that $n< n^{*}, X_n=x >x_n$ and $X_n$ is a local maximum. Then the strategy $(F,S)$ is the pure Nash equilibrium in the game described above.
\end{lemma}

\begin{lemma}
Suppose that the current state of the process $(n,X_n)$ is such that $n< n^{*}, X_n=x \leq x_n$ and $X_n$ is a local maximum. Then the strategy $(F,F)$ is the pure Nash equilibrium in the game described above.
\end{lemma}

\section{Numerical example.}
\subsection{Value of the game.}
\begin{figure}[H]
\includegraphics[width=0.9\textwidth]{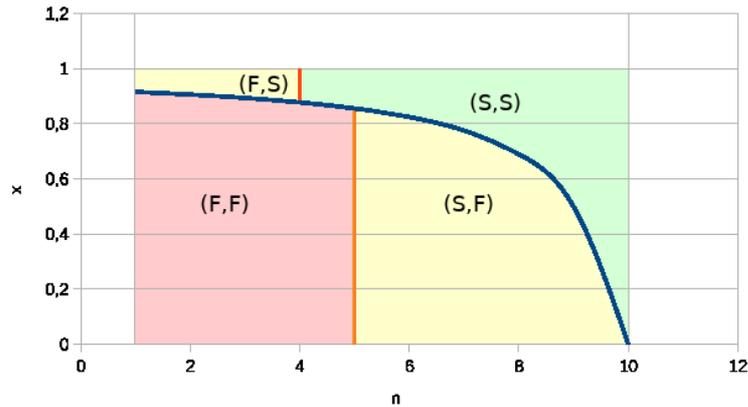}
\caption{Boundaries of the strategies for $N=10, p=0.25$. The shift for the player I is clearly visible. In this case $n^*=4$ but $\tilde{n}=5$. }
\end{figure}
The value of the game for different values of priority parameter $p$ and $N=10$ is presented bellow.
\begin{displaymath}
\begin{split}
(val_{1,10},val_{2,10})=&(-0.00201, 0.19557), \quad p=0.1 \\
(val_{1,10},val_{2,10})=&(0.03283, 0.12896), \quad p=0.25 \\
(val_{1,10},val_{2,10})=&(0.06897, 0.08796), \quad p=e^{-1} \\
(val_{1,10},val_{2,10})=&(0.13662, 0.03787), \quad p=0.5
\end{split}
\end{displaymath}
\subsection{Shift of the threshold for player I}
The table bellow presents different values of the threshold $\tilde{n}$ for different horizons and values of the priority parameter $p$.

\begin{table}[H]
\centering
\begin{tabular}{ |c|c|c|c|c|c|c|c| } 
 \hline
 $N$ & $n^*$ & \multicolumn{6}{c|}{$\tilde{n}$}\\ \hline
 &&$p=0.1$ & $p=0.2$ & $p=0.25$ & $p=1/3$ & $p=e^{-1}$ & $p=0.5$\\
 \hline
 5 & 3 &	3 &	3 &	3 &	3 &	3 &	3\\
10 &	4 &	4 &	5 &	5 &	5 &	5 &	6 \\
20 &	8 &	9 &	10 &	10 &	11 &	11 &	12\\
30 &	12 &	14 &	15 &	15 &	16 &	17 &	18\\
50 &	19 & 24	&	26 & 26	& 28	&	28 &	31\\
 
 \hline
\end{tabular}\caption{Numbers $\tilde{n}$}
\end{table}
\section{Conclusion}
The model presented in this work was created as a fruit of reflection on real problems in the field of business and finance. In the competition between two opponents from which one of them has access to more data we have found the equilibrium states. If the priority parameter of no-information player $p\leq 0.5$ we have found that no-information player has to change his strategy in relation to the situation if he remained in the game alone. However the full-information player does not intend to change his strategy. The numerical examples presented here are good presentation of the model.
\bibliography{mybibfile}
%\bibliographystyle{ws-rv-van}
%\bibliographystyle{siam}
%\nocite{choRobSig71:Expectation}
\bibliographystyle{ws-rv-van}
%\bibliography{ws-rv-sample}
%\printindex[aindx]                % to print author index
\printindex

\end{document}